\documentclass[preprint,12pt]{elsarticle}

\usepackage{amssymb}
\usepackage{rotate}
\usepackage{graphicx}
\usepackage[T1]{fontenc}
\usepackage{fancyhdr}
\usepackage{afterpage}
\usepackage{mathptmx}      
\usepackage{mathrsfs}
\usepackage{amsmath}
\usepackage{latexsym}

\newtheorem{theorem}{\textbf{Theorem}}[section]

\newtheorem{proposition}{\textbf{Proposition}}[section]
\newtheorem{definition}{\textbf{Definition}}[section]
\newdefinition{remark}{Remark}[section]
\newtheorem{example}{\textbf{Example}}[section]

\newproof{proof}{\textbf{Proof}}
\newproof{pot}{Proof of Theorem \ref{thm2}}

\journal{Fuzzy Sets and Systems}

\begin{document}

\begin{frontmatter}



\title{Vague partition, fuzzy sets and their axiomatical foundations}

\renewcommand{\thefootnote}{\fnsymbol{footnote}}
\author{Xiaodong Pan \footnote{Corresponding author.\\
E-mail address: xdpan1@163.com(Xiaodong Pan)}}

\address{School of Mathematics, Southwest Jiaotong University, West Section, High-tech Zone, Chengdu, Sichuan, 611756, P.R. China}

\begin{abstract}
Based on the in-depth analysis of the nature and features of vague phenomenon, this paper focuses on establishing the axiomatical foundation of the membership degree theory for vague phenomenon, presents an axiomatic system of governing membership degrees and their interconnections. Some important basic notions, such as vague variable, vague partition etc. are defined, their useful properties are characterized. Moreover, the notion of fuzzy set is also redefined by the notion of vague partition on the basis of the axiomatic system. Hence, this work can serve as a mathematical model of dealing with the phenomena of vagueness by axiomatical approach from the many-valued point of view, as well as the axiomatical foundation of fuzzy sets and its applications. The thesis defended in this paper is that the difference among vague attributes is the key point to recognize and model vague phenomena, membership degrees should be considered in a vague membership space. In other words, vagueness should be treated from a global and overall (instead of local) of view.
\end{abstract}


\end{frontmatter}


\section{What is vagueness?}
In the 1902 Dictionary of Philosophy and Psychology \cite{Peirce02}, Charles Sander Peirce explained the entry "vague" as follows:
\par \textit{A proposition is vague when there are possible states of things concerning which it is intrinsically uncertain whether, had they been contemplated by the speaker, he would have regarded them as excluded or allowed by the proposition. By intrinsically uncertain we mean not uncertain in consequence of any ignorance of the interpreter, but because the speaker's habits of language were indeterminate. (Peirce 1902, 748)}
\par From the above description, we can see that vagueness is intrinsically indeterminacy, which means such a kind of state of thing that is not absolutely affirmative or negative to say which attribute should be assigned to this state. Vagueness is related closely to the use of natural language, and results in borderline cases. Here it also needs to emphasize that the notion of vagueness in this paper is different from that of Zadeh \cite{Zadeh1978} considered, he wrote "\textit{Although
the terms fuzzy and vague are frequently used interchangeably in the literature, there is, in fact, a significant difference between them. Specifically, a proposition, $p$, is fuzzy if it contains words which are labels of fuzzy sets; and $p$ is vague if it is
both fuzzy and insufficiently specific for a particular purpose.}" In what follows, the terms fuzzy and vague will be used interchangeably unless otherwise stated. The meaning of vagueness will be explained in further detail below.
\par In \cite{Nov99}, Nov\'{a}k said that "\textit{Characterization of the vagueness phenomenon is fundamental for further development of fuzzy logic as well as its applications $\cdots \cdots$ the key role in the roots of fuzzy logic is played by the vagueness phenomenon}". This is very important. Not only fuzzy logic in narrow sense, the characterization of the vagueness phenomenons is fundamental, and plays a key role in almost all theoretical researches and application fields related to fuzzy sets. Meanwhile, the thinking on vagueness has always been an important topic for philosophers, logicians, mathematicians and artificial intelligence experts. However, what exactly the vagueness is? what exactly the fuzzy phenomenons are? How does the vagueness come to be? So far, seemingly it is still not so clear.
\par Until we arrive at answers to such questions, it still needs to give some necessary explanations about vagueness. At present, most of the treatises concerning the essence and characterization of vagueness come from the philosophy community. For more details, please refer to \cite{Dee10, Egr11, Keefe2000, Pab13, Ron11, Smi97}.
\subsection{Sorites Paradoxes}
The word "sorites" derives from the Greek word for heap, the paradox is so named because of its original characterization, attributed to Eubulides of Miletus. Vague predicates are susceptible to Sorites Paradoxes. As we mentioned above, the extension of a vague notion (or predicate) has borderline cases, is unclear. In other words, there are objects which one cannot say with certainty whether belong to a group of objects which are identified with this concept or which exhibit characteristics that have this predicate. Hence, a vague notion is always susceptible to Sorites Paradox.
\par Many people will agree that almost all vague phenomenons can be attributed to the vagueness of predicates. Consider the predicate "\textit{is a heap}", the paradox goes as follows: consider a heap of wheat from which grains are individually removed. One might construct the argument, using premises, as follows:
\begin{itemize}
\item Premise 1: $10^{10}$ grains of wheat is a heap of wheat;
\item Premise 2: when we remove one grain of wheat from a heap of wheat, the remain is still a heap.
\end{itemize}
Repeated applications of Premise 2 (each time starting with one fewer grain) eventually forces one to accept the conclusion that a heap may be composed of just one grain of wheat (and consequently, if one grain of wheat is still a heap, then removing that one grain of wheat to leave no grains at all still leaves a heap of wheat; indeed a negative number of grains must also form a heap). There are some other variations of the paradox, in fact, this paradox can be reconstructed for a variety of predicates, for example, with "\textit{is tall}", "\textit{is rich}", "\textit{is old}", "\textit{is blue}", "\textit{is bald}", and so on.
\par The root of the paradox lies in the existence of borderline cases brought about by the use of vague predicates which are usually expressed by natural language. Even if you may know all the outside information about the vague predicate, such as, the exact number of these wheats, the precise value of the height of a man, the precise value of the wavelength of one color, etc., you are still unable to decide with certainty whether these wheats form a heap, whether this man is a tall man, whether this kind of color is red. That is not because we can't understand these predicates, but because of the vagueness of these predicates.

\subsection{The basic features of vagueness}
Besides Sorites Paradoxes, many people will also agree that almost all vague predicates also share two other features.
\par One is that vague predicates admit borderline cases. Borderline cases are the cases where it is unclear whether the predicate applies. In other words, this means that a vague statement is inconsistent with the law of excluded middle. When we say that some statement satisfies the law of excluded middle, actually, it means that there exists no any vagueness about the statement. It is important to stress that here the borderline refers to the "absolute borderline", that is to say, it can't be got rid of by improved technical methods like observation methods, description methods, language skills etc., these outside methods. For example, the predicate "\textit{is tall}" is vague because a man who is 1.8 meters in height is neither clearly tall nor clearly non-tall. No amount of conceptual analysis or empirical investigation can settle whether a 1.8-meter man is tall, as we mentioned above, it is a kind of intrinsically indeterminacy.
\par The other is that vague predicates apparently lack well-defined extensions. Vagueness arises during the process of grouping together objects according to some property $\varphi$ (expressed by a vague predicate) of objects. The "class" $X$ of all objects which have the property $\varphi$ cannot be taken as a set (crisp set) since the property $\varphi$ makes it impossible for us to characterize the "class" precisely and unambiguously. For example, let the set of all British males be the domain of discourse, consider the "class" of "\textit{tall man}". It is impossible to find a sharp boundary between all tall men and others.
\par It is easy to find out that these three features of vague predicates are also closely related to each other. You could also say that these three features characterize the vagueness of vague predicates from three different perspectives respectively.
\par In addition, the vagueness with respect to (or caused by) a vague predicate $\varphi$ doesn't stem from the shortcomings of the cognitive abilities of the human body, but from a kind of objective attribute exists in the objects itself and human subjective cognitive style. In other words, there is an intermediary transition state between the differences, even very small differences in $\varphi$, the existence of this state makes it impossible for us to present a clearly partition of objects based on the predicate $\varphi$. There exists still an intermediary transition between 1.8 meters and 1.801 meters, one cannot change his (or her) height from 1.8 meters to 1.801 meters abruptly, it should be a continuous changing process.
\subsection{Multidimensional vagueness}
So far we have considered only these vague predicates which are determined by a single dimension of variation (or only one attribute), such as height for "\textit{tall}", age for "\textit{young}" and temperature for "\textit{hot}". But in practical issues, many vague predicates are multidimensional: several different dimensions of variation (or several attributes) are involved in determining their applicability. The applicability of "\textit{big}", when it is used to describe a man, depends on both height and volume. Whether a ball is counted as a "\textit{small red ball}" depends not only on the volume of ball but also on its color. Moreover, there are still some vague predicates in which it is not even a clear-cut set of dimensions determining the applicability of these predicates: it is not clear which factors are related to each other and blend into one another. For example, the applicability of "\textit{good}", which is used to describe undergraduates, whether a undergraduate should be counted as a "\textit{good}" undergraduate, maybe different people think quite differently on it, many factors need to be considered, and it is very difficult to present a universally accepted criterion to judge whether a student is good. Of course, multidimensional vague predicates also share all these features mentioned above of vague predicates.
\par It is important to note, however, that when one applies fuzzy sets to model multidimensional vagueness, then the range of membership function should be a subset in $[0,1]^{n}$ instead of $[0,1]$, where $n$ denotes the number of dimensions of a multidimensional vague predicate. More generally, one can take a complete distributive residuated lattice \cite{War38} as the range of membership function. More details please refer to \cite{Gog67, Gog69}.
\subsection{Higher-order vagueness}
Let $F$ be a vague predicate, there is no (sharp) boundary between the objects that determinately satisfy $F$ and those that do not determinately satisfy $F$, that is, it admits borderline case where it is unclear whether $F$ applies. This is the so-called first-order vagueness. Intuitively, it seems that there is also no (sharp) boundary between the objects that determinately satisfy $F$ and the borderline objects in the borderline case, nor is there a sharp boundary between the borderline objects and those that do not determinately satisfy $F$. This is called as second-order vagueness. Proceeding in this way, one can define the notions of "third-order vagueness", "fourth-order vagueness" etc., which are uniformly referred to as "higher-order vagueness".
\par When one applies fuzzy sets to model higher-order vagueness, the problem is how to define the membership functions of fuzzy sets. For example, let $\varphi$ be the statement "\textit{Tom is tall, if Tom is 1.8 metres in height}", how much is the truth value of $\varphi$? If one set its truth value equal to 0.6 (or other any definite value), then you might ask: "why not 0.61 or 0.59?" Meanwhile, another question also arises: if we consider that the truth value of "\textit{John is tall}" is bigger than 0.6, then how much John's height should be? This is the reflection of higher-order vagueness in fuzzy set theory.
\par In light of fuzzy set theory, several solutions have been proposed by Zadeh \cite{Zad75-1}, Grattan-Guiness \cite{Gra75} and others \cite{Gor87, Hu14, HuW14, Kar01, Men07, Sad13, Zadeh1975} in response to the above problems. Zadeh took linguistic terms as the grades of membership of fuzzy sets, and then these linguistic terms were modeled by fuzzy sets of type $n$ whose membership function ranges over fuzzy sets of type $n - 1$. The membership function of a fuzzy set of type 1 ranges over the interval $[0, 1]$. Another way to deal with the higher vagueness in fuzzy set theory is to replaced the interval $[0, 1]$, the set of grades of membership of fuzzy sets, by the set of subintervals of $[0, 1]$, which has been proposed by Grattan-Guiness. I. in 1975. However, it should also be pointed out that the above two strategies can only alleviated the problem of higher vagueness to some extent, but it can not solve the problem completely. In fact, for fuzzy sets of type $n$, we still need to define the membership function of fuzzy sets of type $n - 1$. For interval valued fuzzy sets, we still need to define the subintervals which are the grades of membership of fuzzy sets. But how to do that? And how to estimate it's rationality? It's still not very clear.
\subsection{Vagueness vs. Ambiguity, Relativity and Underspecificity}
Next, we present some objects which are easily mistaken for vagueness. First of all, ambiguity is not vagueness. Certainly,
One term can be ambiguous and vague: "bank" for example has two quite different main senses (concerning financial institutions or river
edges), both of which are also vague. One term can also be vague but not ambiguous. Look at the following example, it is natural to suppose that the predicate "is tall man" has a univocal sense (is not ambiguous),  but this predicate does not determine a sharp, well-defined extension. Usually, you can't answer questions such as "where is the bank", even if you ignore the vagueness of the term "bank". That only mean that you don't know what the term "bank" in the question refers to, there exists ambiguity which is due to the absence of the context about the question. Meanwhile, you can't also answer questions such as "is Tom (his height is 1.8 meters) a tall man". Here the case is quite different, you can't answer is only due to the vagueness in the predicate "\textit{is tall}".
\par Secondly, vagueness is different from relativity. Some researchers \cite{Sha06} think that vagueness should be straightforwardly identified with paradigm context dependence (i.e. having different extensions in different contexts), even though many terms (e.g. "\textit{is tall}") have both of the features, vagueness and relativity. However, when one fixes the context which can be made as definite as he (or she) likes (in particular, choose a specific comparison class, e.g. current British professional soccer players), the predicate remains vague, admits borderline cases and has fuzzy boundaries, and Sorites Paradox will retain its force. This indicates that we are unlikely to understand vagueness or solve the paradox by concentrating on context-dependence. In the following, unless otherwise stated, the discussion about vagueness is always carried out in a specific context.
\par Thirdly, we don't also consider expressions like "someone said something" as a vague expression. Its meaning is not specific just only because we don't have adequately information for purposes in hand, which is also called underspecificity. Other examples such as "$X$ is an integer smaller than one thousand" are also in this case. The vague predicates mentioned in this paper will always be in a specific context, and its meaning are transpicuous and understandable.
\subsection{Vagueness vs. Uncertainty, Randomness}
Vagueness is also often confused with uncertainty and randomness. In general, uncertainty is applied to predictions of future events, to physical measurements that are already made, or to the unknown. It arises in partially observable and/or stochastic environments, as well as due to ignorance and/or indolence. Douglas W. Hubbard \cite{Hubbard2010} defined uncertainty as: \textit{"The lack of complete certainty, that is, the existence of more than one possibility. The "true" outcome/state/result/value is not known"}. Uncertainty is encountered when an experiment (process, test, etc.) is to proceed, the result of which is not known to us. Hence, uncertainty is always connected with the question whether the given event may be regarded within some time period, or not; there is no uncertainty after the experiment was realized and the result is known to us.
\par Uncertainty emerges probably due to the lack of enough knowledge, probably due to the shortcomings of our cognitive abilities, and also probably due to the relatively poor technical conditions etc. Frequently, uncertainty will disappear as long as these situations have been improved. Here you see that uncertainty differs from vagueness, the latter has nothing to do with these outside conditions. Vagueness only concerns the object itself under consideration and the way how it is delineated according to its certain attribute, and won't disappear as time passed. You can say that the difference between uncertainty and vagueness corresponds to the difference between potentiality and factuality.
\par In addition, Randomness is a specific form of uncertainty, and thus differs from vagueness. "Randomness" means according to Wikipedia "\textit{Having no definite aim or purpose; not sent or guided in a particular direction; made, done, occurring, etc., without method or conscious choice; haphazard.}" This concept suggests a non-order or non-coherence in a sequence of symbols or steps, such that there is no intelligible pattern or combination. Probability theory is the mathematical analysis of random phenomena, probability can be thought of as a numerical measure of the likelihood that an event will occur, its value is a number between 0 (0 percent chance or will not happen) and 1 (100 percent chance or will happen). In this sense, probability is very similar to the membership degrees in fuzzy set theory proposed by Zadeh \cite{Zadeh1965} in 1965, which has been viewed a suitable tool to treat vagueness.
\subsection{What is vagueness?}
Now, let's go back to the previous questions: what exactly the vagueness is? how does it come to be? Many people think that vagueness stems from the use of natural language, is this really the case? Are all kinds of vaguenesses linguistic?
\par Based on the previous discussion, our point of view is that vagueness arises in the process of classifying objects, it is a kind of manifestation of the continuity and gradualness existing in the process of development and evolution of objects. As an example, one person is impossible to become an adult from a baby in an instant, the process is a continual and evolutionary. When one regards this process in a discrete point of view based on some specific purpose, or divides this process into several discrete stages, consequently the continuity and gradualness in objects themselves can not be revealed completely under such partition way, and then the vagueness arises.
\par An usual way to discretize the continual and evolutionary process of objects is to delineate the process by using some natural linguistic units (words or phrases), each natural linguistic unit represents a stage in the process, they summarize the main characteristic of the corresponding stage respectively. In other words, natural language provides us with a very useful tool to describe the continuous process in a discrete way. Hence, vagueness usually arises together with the use of natural language.
\par Consider the following example, let the set of all Chinese be the domain of discourse, and consider the attribute "\textit{Age}". According to our common sense, one person's age should be in the scope of 0 to 200 (an optimistic estimate). It is a continuous process for one's age to change from 0 to 100, if you divide the process into several discrete stages, and label each stage with one of the following natural linguistic units e.g. \textit{Infant, Childish, Juvenile, Youth, Adult, Middle age, Elderly, Old, Senile} etc., then this kind of discrete partition hides the inherent continuity and gradualness existing in the attribute "Age". Consequently, these predicates such as "\textit{is young}", "\textit{is adult}", "\textit{is old}" etc., show the vagueness.
\par According to its birth, vagueness can be fallen into different types: one dimensional vagueness, two dimensional vagueness, etc.
\section{The mathematical analysis of vagueness}
In order to deal with various kinds of problems concerning vagueness conveniently and effectively, naturally we hope to model or describe vagueness by a numerical (or formalized) way. One of the most important contributions in this direction is the establishment and development of fuzzy set and fuzzy logic by Zadeh and others \cite{Bed13, Cin11, Dubois99, Haj98, Nov99, Nov06, Nov12, Pav79, Zadeh1965}. Fuzzy sets are designed to model the extensions of vague concepts, such as \textit{Youth}. In other words, the extension of a vague concept $\varphi$ is taken to be a fuzzy set $A$, which is defined by a membership function $\mu_{\varphi}: U \rightarrow [0,1]$ introduced by Zadeh in 1965, where $U$ is the domain of discourse. For each element $x \in U$, the value $\mu_{\varphi}(x) \in [0,1]$ is the membership degree of $x$ in $A$, or the degree of truth that the element $x$ belongs to $A$. Membership functions is a kind of generalization of the characteristic functions of ordinary sets. The crucial difference between ordinary sets and fuzzy sets is thus in using of the scale $[0, 1]$.
\par As mentioned in the Section 1, vagueness is a kind of manifestation of the continuity and gradualness in the process of development and evolution of objects. The graded approach developed in fuzzy sets and fuzzy logic reproduces, at least to some extent, this kind of continuity and gradualness, then the objects are associated with numbers in such a way that the relation between objects and its property is faithfully represented as numerical properties. For example, let $\varphi$ be the property "\textit{tall}" and $x$ a concrete object, the number $\mu_{\varphi}(x)$ captures the degree of intensity which $x$ possesses the property $\varphi$.
\par However, is fuzzy set or fuzzy logic the suitable and effective tool for dealing with vagueness? Seemingly the answer is not so obviously. In \cite{Beh09}, B\v{e}hounek said that "\textit{Fuzzy logic cannot claim to be the logic of vagueness, most of which are not captured by deductive fuzzy logic}". Nov\'{a}k said in \cite{Nov12} that "\textit{fuzzy logic is not the logic of vagueness but the logic of ordered structure $\cdots$ there are not yet many results on the applied side $\cdots$ It is still not fully clear which logic is the most convenient to solve problems related to models of vagueness and their applications}". You can see from these comments that fuzzy set theory is still at its first stage, its theoretical foundation is far from complete. There is still many works that need to be done, some important foundational notions and operations in fuzzy set theory and fuzzy logic still require further explanation, and still need further to be more normalized, explicitly and strictly. For example, how to define the membership function of a given fuzzy set? how to define the operations between fuzzy sets? and how to estimate the rationality of various application methods based on fuzzy sets? And so on.
\par This paper aims at doing some attempts in this direction. Based on these studies about vagueness in Section 1, we think that for the characterizing of a vague predicate, we should use the relationship between this predicate and other closely related vague predicates as the starting point. As we all know, cognition comes from comparison. If you want to know what red is, then you need also to know what yellow, orange and other colors are. If you want to know what a tall man is like, you first need to know what a short man is like. Differences are the fundamental of understanding things. From a global and overall of view, we will try to establish a mathematical model of treating phenomenons of vagueness by axiomatical approach, it is our hope that the model could serve as the theoretical foundation of fuzzy set theory and its application, and further serve as the starting point of formalized theory of dealing with the phenomenons of vagueness.
\section{Some basic terminologies}
In this section, we define the objects of our future study: elementary vague attribute, vague attribute, vague space and vague judgement.
\par In general, a concept (of course, it is crisp or distinct) is determined by or involved with one or a group of attributes, these attributes can be one-dimensional, and can also be multidimensional. For example, \textit{Weight, Height, Age} are one-dimensional attributes with respect to the concept "Man", \textit{Length, Width, Height, Color} are also one-dimensional attributes with respect to the concept "Box", but \textit{Area, Volume} are multidimensional attributes with respect to the concept "Room". It is worth noting that a multidimensional attribute can be determined usually by several one-dimensional attributes. For example, $Area$ is determined by one-dimensional attributes \textit{Length} and \textit{Width}, and can be regarded as the Cartesian product $Length \times Width$.
\par Given a concept $C$, consider the partition (finite or countable) of its extension (which is usual a crisp (or classical) set or class) according to one or a group of attributes with respect to $C$, the partition is not sharply if the attributes as the standard of the partition is characterized by continuity, and then vagueness arises, as mentioned in Section 1.7. For example, consider the partition of the extension of the concept "Man" according to the attribute "\textit{Height}", it can be divided into the following "classes": \textit{Dwarf, Short, Medium, Tall,} etc.
\par As the first step, we consider only the partition according to certain one-dimensional attribute. In general, when we speak of the height of somebody, we might prefer to say that he (or she) is \textit{tall}, \textit{short}, or \textit{medium} etc., rather than to say that he (or she) is 1.753 metres in height. Let $C$ be the concept "Man" and $\varphi_{C}$ the one-dimensional attribute "\textit{Height}". The set $$\{\text{Dwarf, Very Short, Short, Medium, Tall, Very Tall, Very Very Tall}\}$$ is called the elementary vague attribute set of $\varphi_{C}$, denoted by $\Omega_{\varphi_{C}}$. In fact, the set $\Omega_{\varphi_{C}}$ can be viewed as a vague or fuzzy partition of the extension of $C$ according to the attribute $\varphi_{C}$. Here it is worth noting that the set $\Omega_{\varphi_{C}}$ is changeable depends on the need of the partition in practical problems. The formalized definition of vague partition will be presented in the latter section. In what follows, we always denote the elementary vague attribute set of $\varphi_{C}$ with respect to the concept $C$ by $\Omega_{\varphi_{C}}$ (which can be effectively equivalent to the sample space in probability theory, the only difference is that in vague membership theory, the elementary vague attribute set is alterable according to different requirements of applications), it can be finite or countable infinite, even uncountable. The elements in $\Omega_{\varphi_{C}}$ are called elementary vague attributes with respect to $\varphi_{C}$. In addition, the domain of discourse of the attribute $\varphi_{C}$, or the range of values of measurement of $\varphi_{C}$, is denoted by $U_{\varphi_{C}}$. Based on the analysis in Section 1, the set $U_{\varphi_{C}}$ should be continuous. For example, let $C$ be the concept "Man" and $\varphi_{C}$ the attribute "\textit{Age}", then you can take the universe $U_{\varphi_{C}}$ to be the interval $[0, 300]$ according to the common sense.
\par Let $\mathcal{F} = \{\bot, \top, \neg, \barwedge, \veebar\}$, where $\bot, \top$ represent the nonexistent vague attribute and the intrinsic vague attribute respectively, and $\neg, \barwedge$ and $\veebar$ represent the connectives "negation", "and" and "or" respectively, which are used as the operations among various vague attributes. The vague attribute set of the concept $C$ with respect to the attribute $\varphi_{C}$, denoted by $\Sigma_{\varphi_{C}}$, is the smallest set such that
\begin{itemize}
\item $\bot, \top \in \Sigma_{\varphi_{C}}$.
\item $\Omega_{\psi_{C}} \subset \Sigma_{\varphi_{C}}$.
\item If $A \in \Sigma_{\varphi_{C}}$, then $\neg A \in \Sigma_{\varphi_{C}}$.
\item If $A, B \in \Sigma_{\varphi_{C}}$, then $A \barwedge B, A \veebar B \in \Sigma_{\varphi_{C}}$.
\item If $A_{1}, A_{2}, \cdots, A_{n}, \cdots \in \Sigma_{\varphi_{C}}$, then $\veebar_{i = 1}^{\infty}A_{i}, \barwedge_{i = 1}^{\infty}A_{i} \in \Sigma_{\varphi_{C}}$.
\end{itemize}
In what follows, we call the triple $(\Omega_{\varphi_{C}}, \Sigma_{\varphi_{C}}, \mathcal{F})$ vague space with respect to $\varphi_{C}$ (which can be effectively equivalent to the event field in probability theory), and call the elements in $\Sigma_{\varphi_{C}}$ vague attributes with respect to $\varphi_{C}$. You may have noticed that a vague attribute can be in some sense regarded as an evaluating linguistic expression which has been introduced by Dvo\v{r}\'{a}k and Nov\'{a}k in \cite{Dvorak04}. As mentioned above, $Short$ is an elementary vague attribute with respect to "\textit{Height}", the element $Short \veebar Medium$ is also a vague attribute with respect to "\textit{Height}", which means "\textit{Short or Medium}", and can be named as the attribute "\textit{Lower Medium}".
\par For any $x \in U_{\varphi_{C}}$, a vague judgement in vague space $(\Omega_{\varphi_{C}}, \Sigma_{\varphi_{C}}, \mathcal{F})$ with respect to $x$ is a procedure to determine these degrees to which these vague attributes in $\Omega_{\varphi_{C}}$ are possessed by the element $x$ respectively. Put it in another way, a vague judgement in vague space $(\Omega_{\varphi_{C}}, \Sigma_{\varphi_{C}}, \mathcal{F})$ with respect to $x$ is a series of statements (vague propositions) "$x$ is $p$", where $p \in \Omega_{\varphi_{C}}$. For instance, let $C$ be the concept "Man" and $\varphi_{C}$ the attribute "\textit{Height}". Define $\Omega_{\varphi_{C}}$ and $U_{\varphi_{C}}$ as $$\{\text{Dwarf, Very Short, Short, Medium, Tall, Very Tall, Very Very Tall}\}$$ and $[0, 3]$, respectively, and $\mathcal{F} = \{\bot, \top, \neg, \barwedge, \veebar\}$ and $\Sigma_{\varphi_{C}}$ is given accordingly. If $x = 1.7 \ \text{metres} \in U_{\varphi_{C}}$, then a vague judgement in vague space $(\Omega_{\varphi_{C}}, \Sigma_{\varphi_{C}}, \mathcal{F})$ with respect to $x$ consists of the following statements: "$x$ \textit{is dwarf}"; "$x$ \textit{is very short}"; "$x$ \textit{is short}"; "$x$\textit{ is medium}"; "$x$\textit{ is tall}"; "$x$ \textit{is very tall}"; "$x$\textit{ is very very tall}". Furthermore, we can assign a real number in the unit interval $[0, 1]$ to each of these statements, which are referred to as the truth values of these statements. The truth value of the statement "$x$ is $p$" ($p \in \Omega_{\varphi_{C}}$) is a measure or estimation of the extent which one trusts that $x$ is $p$. According to fuzzy set theory, this truth value is also called membership degree of $x$ in fuzzy set $A$, where $A$ is determined by the vague attribute $p$, and is a mathematical model of the extension of $p$. In a vague judgement, the relationship among membership degrees of vague propositions is the key point this paper focuses on, we will discuss this topic in detail in the next section.
\section{The axioms for membership degrees}
In this section, we present the system of axioms that will govern the relations among various degrees of membership. As a preliminary step, we first review several fundamental concepts, "Triangular norm", "Triangular conorm" and "Strong negation", which will be needed for the rest part of this paper. For more details, please refer to \cite{Klement}.
\begin{definition}
A triangular norm (t-norm for short) is a binary operation $T$ on the unit interval $[0,1]$, i.e., a function $T: [0, 1]^{2} \rightarrow [0, 1]$, such that for all $x, y, z \in [0, 1]$, the following four axioms are satisfied:\\
(T1) \ $T(x, y) = T(y, x)$,\hspace*{\fill} (commutativity)\\
(T2) \ $T(x, T(y, z)) = T(T(x, y), z)$,\hspace*{\fill} (associativity)\\
(T3) \ $T(x, y) \leqslant T(x, z)$ whenever $y \leqslant z$,\hspace*{\fill} (monotonicity)\\
(T4) \ $T(x, 1) = x$.\hspace*{\fill} (boundary condition)
\end{definition}
\begin{example}
The following are the four basic t-norms $T_{M}, T_{P}, T_{L},$ and $T_{D}$ given by, respectively:\\
\par $T_{M}(x, y) = \min\{x, y\}$,\hspace*{\fill} (minimum)
\par $T_{P}(x, y) = x \cdot y$,\hspace*{\fill} (product)
\par $T_{L}(x, y) = \max\{x + y - 1, 0\}$,\hspace*{\fill} ({\L}ukasiewicz t-norm)
\par $T_{D}(x, y) = \left\{ \begin{array}{ll}
0, &\mbox{if \ $(x, y) \in [0, 1)^{2}$,} \\
\min\{x, y\}, &\mbox{otherwise.}
\end{array} \right. $  \hspace*{\fill} (drastic product)
\end{example}
\begin{definition}
A triangular conorm (t-conorm for short) is a binary operation $S$ on the unit interval $[0,1]$, i.e., a function $S: [0, 1]^{2} \rightarrow [0, 1]$, which, for all $x, y, z \in [0, 1]$, satisfies ($T1$)-($T3$) and \\
(S4) \ $S(x, 0) = x$. \hspace*{\fill} (boundary condition)
\end{definition}
\begin{example}
The following are the four basic t-conorms $S_{M}, S_{P}, S_{L},$ and $S_{D}$ given by, respectively:\\
\par $S_{M}(x, y) = \max\{x, y\}$,\hspace*{\fill} (maximum)
\par $S_{P}(x, y) = x + y - x \cdot y$,\hspace*{\fill} (probabilistic sum)
\par $S_{L}(x, y) = \min\{x + y, 1\}$,\hspace*{\fill} ({\L}ukasiewicz t-conorm, bounded sum)
\par $S_{D}(x, y) = \left\{ \begin{array}{ll}
1, &\mbox{if \ $(x, y) \in (0, 1]^{2}$,} \\
\max\{x, y\}, &\mbox{otherwise.}
\end{array} \right. $  \hspace*{\fill} (drastic sum)
\end{example}
\par Since both t-norm and t-conorm are algebraic operations on the unit interval $[0, 1]$, it is of course also acceptable to use infix notations like $x \otimes y$ and $x \oplus y$ instead of the prefix notations $T(x, y)$ and $S(x, y)$ respectively. In what follows, we will use these infix notations most of the time.
\begin{definition}
(i) A non-increasing function $N: [0, 1] \rightarrow [0, 1]$ is called a negation if
 \par (N1) $N(0) = 1$ and $N(1) = 0$.\\
(ii) A negation $N: [0, 1] \rightarrow [0, 1]$ is called a strict negation if, additionally,
 \par (N2) $N$ is continuous.
 \par (N3) $N$ is strictly decreasing.\\
(iii) A strict negation $N: [0, 1] \rightarrow [0, 1]$ is called a strong negation if it is an involution, i.e., if
 \par (N4) $N \circ N = id_{[0,1]}$ is continuous.
\end{definition}
\par It is obvious that $N: [0, 1] \rightarrow [0, 1]$ is a strict negation if and only if it is a strictly decreasing bijection.
\begin{example}
(i) The most important and most widely used strong negation is the standard negation $N_{S}: [0, 1] \rightarrow [0, 1]$ given by $N_{S} = 1 - x$. It can be proved that each strong negation can be seen as a transformation of the standard negation, see \cite{Klement}.\\
(ii) The negation $N: [0, 1] \rightarrow [0, 1]$ given by $N(x) = 1 - x^{2}$ is strict, but not strong.\\
(iii) An example of a negation which is not strict and, subsequently, not strong, is the G$\ddot{o}$del negation $N_{G}: [0, 1] \rightarrow [0, 1]$ given by
\[
N_{G}(x) = \left\{
\begin{array}{ll}
1, &\mbox{if \ $x = 0$,}\\
0, &\mbox{if \ $x \in (0, 1]$.}
\end{array}
\right.
\]
\end{example}
\begin{definition}
If $T$ is a t-norm, then its dual t-conorm $S: [0, 1]^{2} \rightarrow [0, 1]$ is given by $S(x, y) = 1 - T(1 - x, 1 - y)$.
\end{definition}
\par In the following, we list some basic properties of t-norm, t-conorm and negation which will be needed in the sequel.
\begin{proposition}
Let $T$ be a t-norm. Then the following conclusions hold:\\
(1) $T(0, x) = T(x, 0) = 0$, T(1, x) = x.\\
(2) $T(x_{1}, y_{1}) \leqslant T(x_{2}, y_{2})$ whenever $x_{1} \leqslant x_{2}$ and $y_{1} \leqslant y_{2}$.\\
(3) The Minimum $T_{M}$ is the only t-norm satisfying $T(x, x) = x$ for all $x \in (0, 1)$.\\
(4) The drastic product $T_{D}$ is the only t-norm satisfying $T(x, x) = 0$ for all $x \in (0, 1)$.
\end{proposition}
\par For t-conorms, the corresponding conclusions can also be obtained similarly.
\begin{definition}
(i) If, for two t-norms $T_{1}$ and $T_{2}$, the inequality $T_{1}(x, y) \leqslant T_{2}(x, y)$ holds for all $(x, y) \in (0, 1)^{2}$, then we say that $T_{1}$ is weaker than $T_{2}$ or, equivalently, that $T_{2}$ is stronger than $T_{1}$, and we write in this case $T_{1} \leqslant T_{2}$.\\
(ii) We shall write $T_{1} < T_{2}$ whenever $T_{1} \leqslant T_{2}$ and $T_{1} \neq T_{2}$, i.e., if $T_{1} \leqslant T_{2}$, but $T_{1}(x_{0}, y_{0}) < T_{2}(x_{0}, y_{0})$ holds for some $(x_{0}, y_{0}) \in [0, 1]^{2}$.
\end{definition}
\begin{proposition}
(1) $T_{D} \leqslant T \leqslant T_{M}$ for arbitrary t-norm $T$.
(2) $T_{D} < T_{L} < T_{P} < T_{M}$.
\end{proposition}
\begin{remark}
The above partial order, which is similar to that in Definition 4.5, can also be defined over all the t-conorms. By the duality between t-norms and t-conorms and the above proposition, we have that $S_{D} \geqslant S \geqslant S_{M}$ for arbitrary t-conorm $S$ and $S_{D} > S_{L} > S_{P} > S_{M}$.
\end{remark}
\par In order to deal with phenomenons of vagueness mathematically (numerically or formalized), and establish the rigorous foundation for the mathematical analysis of vagueness, we propose the following axioms by which the membership degrees and their interconnections are to be governed.
\par First of all, we only consider one dimensional vague attributes. Let $C$ be a concept, $x \in U_{\varphi_{C}}$ and $(\Omega_{\varphi_{C}}, \Sigma_{\varphi_{C}}, \mathcal{F})$ a vague space with respect to the attribute $\varphi_{C}$. The vague membership space with respect to $\varphi_{C}$ (which is one dimensional) and associated with $x$ is a quintuple $(\Omega_{\varphi_{C}}, \Sigma_{\varphi_{C}}, \mathcal{M}_{x}, \mathcal{F}, \mathcal{T})$, where $\mathcal{T} = \{N, \oplus, \otimes\}$ consists of a t-norm $\otimes$, a t-conorm $\oplus$ and a strong negation $N$. $\mathcal{M}_{x}$ is a real value function (it is usually called vague membership measure on vague space $(\Omega_{\varphi_{C}}, \Sigma_{\varphi_{C}}, \mathcal{F})$, or membership measure for short) from $\Sigma_{\varphi_{C}}$ to $[0, 1]$ with respect to $x$, satisfies the following axioms:
\begin{itemize}
\item {\bf Axiom I.} For any $A \in \Sigma_{\varphi_{C}}$, $0 \leqslant \mathcal{M}_{x}(A) \leqslant 1$, and there is at least one element $p \in \Omega_{\varphi_{C}}$ such that $\mathcal{M}_{x}(p) > 0$. If there is $p_{0} \in \Omega_{\varphi_{C}}$ such that $\mathcal{M}_{x}(p_{0}) = 1$, then for any $p \in \Omega_{\varphi_{C}}, p \neq p_{0}$, satisfies $\mathcal{M}_{x}(p) = 0$.\\
  This means that for any $A \in \Sigma_{\varphi_{C}}$, the degree to which $x$ has the vague attribute $A$ is between 0 and 1, and $x$ has at least one of elementary vague attributes to some extent that is bigger than 0, and $x$ has at most one of elementary vague attribute $p$ over the level of membership degree 1.
\item {\bf Axiom II.} $\mathcal{M}_{x}(\bot) = 0, \mathcal{M}_{x}(\top) = 1$. \\
This means that the degree to which $x$ has the nonexistent vague attribute is 0, and that the degree to which $x$ has the intrinsic vague attribute is 1.
\item {\bf Axiom III.} For any vague attribute $A \in \Sigma_{\varphi_{C}}$, $\mathcal{M}_{x}(\neg A) \leqslant \big(\mathcal{M}_{x}(A)\big)^{N}$, where $N$ is a strong negation on $[0, 1]$.\\
This means that the degree to which $x$ doesn't has the vague attribute $A$ is less than or equal to the result of the strong negation operation on the degree to which it has the vague attribute $A$.
\item {\bf Axiom IV.} The countable sum and countable product: for any vague attributes sequence: $A_{1}, A_{2}, \cdots, A_{n}, \cdots \in \Sigma_{\varphi_{C}}$, $$\mathcal{M}_{x}\big(\veebar_{n = 1}^{\infty}A_{n}\big) = \mathcal{M}_{x}(A_{1}) \oplus \mathcal{M}_{x}(A_{2}) \oplus \cdots \oplus \mathcal{M}_{x}(A_{n}) \oplus \cdots,$$ and $$\mathcal{M}_{x}\big(\barwedge_{n = 1}^{\infty}A_{n}\big) = \mathcal{M}_{x}(A_{1}) \otimes \mathcal{M}_{x}(A_{2}) \otimes \cdots \otimes \mathcal{M}_{x}(A_{n}) \otimes \cdots.$$
    Where $\oplus, \otimes$ are triangular conorm and triangular norm respectively, and they are mutually $N$-dual to each other, and $N$ is same with that of Axiom III. \\
    This means that the degree to which $x$ has the vague attribute $\veebar_{n = 1}^{\infty}A_{n}$ equals the result of the t-conorm operation on all these degrees to which it has each vague attribute in the vague attributes sequence, and that the degree to which $x$ has the vague attribute $\barwedge_{n = 1}^{\infty}A_{n}$ equals the result of the t-norm operation on all these degrees to which it has each vague attribute in the vague attributes sequence.
\item {\bf Axiom V.} For any $p \in \Omega_{\varphi_{C}}$, $\mathcal{M}_{x}(\neg p) \geqslant \bigoplus_{q \in \Omega_{\varphi_{C}} \setminus \{p\}}\mathcal{M}_{x}(q)$. \\
This means that the elements of $\Omega_{\varphi_{C}}$ are mutually exclusive to some extent. Hence, $\Omega_{\varphi_{C}}$ makes a vague partition of $U_{\varphi_{C}}$.
\item {\bf Axiom V'.} $0 < \sum_{p \in \Omega_{\varphi_{C}}}\mathcal{M}_{x}(p) \leqslant 1$.
\end{itemize}
\par The vague membership space $(\Omega_{\varphi_{C}}, \Sigma_{\varphi_{C}}, \mathcal{M}_{x}, \mathcal{F}, \mathcal{T})$ is said to be regular if we take $\mathcal{M}_{x}(\neg A) = \big(\mathcal{M}_{x}(A)\big)^{N}$ in Axiom III. Obviously, here the mapping $\mathcal{M}_{x}$ is determined only by its values on $\Omega_{\varphi_{C}}$. It is easy to show that the system of Axioms I-V are consistent. This can be shown by the following example.
\par Let $U_{\varphi_{C}} = [0, 200]$ and $$\Omega_{\varphi_{C}} = \{[0, 40], (40, 80], (80, 120], (120, 160], (160, 200]\},$$ $\bot = \emptyset$ and $\top = [0, 200]$, $\veebar, \barwedge, \neg$ represent the basic set operations, namely union, intersection and complement with respect to $U_{\varphi_{C}}$, respectively. $\Sigma_{\varphi_{C}}$ is the set field generated by $\Omega_{\varphi_{C}}$. $N$ is the standard negation on $[0, 1]$, $\oplus$ is the maximum or G\"{o}del t-conorm and $\otimes$ is the minimum or G\"{o}del t-norm. Let $x = 25$, define the function $\mathcal{M}_{x}$ as follows: $\mathcal{M}_{x}(\bot) = 0$, $\mathcal{M}_{x}(\top) = 1$, $\mathcal{M}_{x}([0, 40]) = 1$ and $$\mathcal{M}_{x}((40, 80]) = \mathcal{M}_{x}((80, 120]) =\mathcal{M}_{x}((120, 160]) = \mathcal{M}_{x}((160, 200]) = 0.$$ It is easy to show that $\mathcal{M}_{x}$ is a vague membership measure defined on the vague space $(\Omega_{\varphi_{C}}, \Sigma_{\varphi_{C}}, \mathcal{F})$. Hence, $(\Omega_{\varphi_{C}}, \Sigma_{\varphi_{C}}, \mathcal{M}_{x}, \mathcal{F}, \mathcal{T})$ is a regular vague membership space.
\par A vague membership space $(\Omega_{\varphi_{C}}, \Sigma_{\varphi_{C}}, \mathcal{M}_{x}, \mathcal{F}, \mathcal{T})$ is said to be normal if there is an element $p_{0} \in \Omega_{\varphi_{C}}$ such that $\mathcal{M}_{x}(p_{0}) = 1.$  In this case, $x$ is said to be crisp in the vague membership space $(\Omega_{\varphi_{C}}, \Sigma_{\varphi_{C}}, \mathcal{M}_{x}, \mathcal{F}, \mathcal{T})$. In fact, any partition in classical sense of set $X$ is a normal vague membership space with respect to any $x \in X$, as it was shown in the above example. It is easy to show that any normal vague membership space must be regular, not vice versa.
\par When $\bigoplus_{p \in \Omega_{\varphi_{C}}}\{\mathcal{M}_{x}(p)\} = a$, then the number $a$ is called the consistent degree of $x$ with the vague membership space $(\Omega_{\varphi_{C}}, \Sigma_{\varphi_{C}}, \mathcal{M}_{x}, \mathcal{F}, \mathcal{T})$. And let $$b = \min_{x \in U_{\varphi_{C}}}\big\{\bigoplus_{p \in \Omega_{\varphi_{C}}}\{\mathcal{M}_{x}(p)\}\big\},$$ then the number $1 - b$ is called the degree of separation of the set $\Omega_{\varphi_{C}}$ of elementary vague attributes.
\par For any $A, B \in \Sigma_{\varphi_{C}}$, $A$ and $B$ are said to be incompatible in the vague space $(\Omega_{\varphi_{C}}, \Sigma_{\varphi_{C}}, \mathcal{F})$ if for any element $x \in U_{\varphi_{C}}$, $\mathcal{M}_{x}(A \barwedge B) = 0$ in vague membership space $(\Omega_{\varphi_{C}}, \Sigma_{\varphi_{C}}, \mathcal{M}_{x}, \mathcal{F}, \mathcal{T})$; $A$ and $B$ are said to be absolutely incompatible in the vague space $(\Omega_{\varphi_{C}}, \Sigma_{\varphi_{C}}, \mathcal{F})$ if for any element $x \in U_{\varphi_{C}}$, $\mathcal{M}_{x}(A)$ and $\mathcal{M}_{x}(B)$ can't be greater than zero at the same time in vague membership space $(\Omega_{\varphi_{C}}, \Sigma_{\varphi_{C}}, \mathcal{M}_{x}, \mathcal{F}, \mathcal{T})$.
\par For example, let $C$ be the concept "Man" and $\varphi_{C}$ the attribute "\textit{Age}", and $\Omega_{\varphi_{C}}$ is defined as the set $$\{Infant, \ Childish, \ Juvenile, \ Youth, Adult, \ Middle\ age, Elderly, \ Old, \ Senile\}.$$ Then the vague attributes "\textit{Childish}" and "\textit{Middle age}" are absolutely incompatible, this is obvious because that for anyone, he (or she) can't has the attributes "Childish" and "Middle age" at the same time to some extent.
\par The mathematical analysis of phenomenons of vagueness, as a mathematical discipline, can and should be developed from axioms in exactly the same way as Probability, and Geometry and Algebra. This means that after we have defined the elements to be studied and their basic relations, and have stated the axioms by which these relations are to be governed, all further exposition must based exclusively on these axioms, independent of the usual concrete meaning of these elements and their relations.
\par In accordance with the above discussion, in Section 3 the concept of the vague attribute set is defined as a free algebra on the elementary vague attribute set. What the elements of this set represent is of no importance in the purely mathematical development of the theory of vague membership degree.
\par Every axiomatic (abstract) theory admits, as is well known, of an unlimited number of concrete interpretations besides those from which it was derived. Thus we find applications in fields of science which have no relation to the concepts of vague attribute and of membership degree in the precise meaning of these words.
\par In what follows, for the sake of convenience, we denote a vague membership space $(\Omega_{\varphi_{C}}, \Sigma_{\varphi_{C}}, \mathcal{M}_{x}, \mathcal{F}, \mathcal{T})$, which is with respect to the attribute $\varphi_{C}$ of the concept $C$ and associated with $x$, by $(\Omega, \Sigma, \mathcal{M}, \mathcal{F}, \mathcal{T})$ when it does not involve any concrete problems.
\section{Immediate corollaries of the axioms, independence and conditional membership degree}
In the sequel, unless otherwise stated, $(\Omega, \Sigma, \mathcal{M}, \mathcal{F}, \mathcal{T})$ always denotes a vague membership space.
\begin{proposition}(Finite sum and finite product)Let $A_{1}, A_{2}, \cdots, A_{n} \in \Sigma$. Then $$\mathcal{M}\big(\veebar_{i = 1}^{n}A_{i}\big) = \mathcal{M}(A_{1}) \oplus \mathcal{M}(A_{2}) \cdots \oplus \mathcal{M}(A_{n}),$$ and $$\mathcal{M}\big(\barwedge_{i = 1}^{n}A_{i}\big) = \mathcal{M}(A_{1}) \otimes \mathcal{M}(A_{2}) \cdots \otimes \mathcal{M}(A_{n}).$$
\end{proposition}
\begin{proof}
Let $A_{n + 1} = A_{n + 2} = \cdots = \bot $ and $A_{n + 1} = A_{n + 2} = \cdots = \top$ respectively, then the above two equalities are the immediate corollaries of Axiom II and IV by combining Definition 4.1 and 4.2.
\end{proof}
\par By Axioms II-IV, Proposition 5.1, Definition 4.1 and 4.2, it is easy to prove the following equalities.
\begin{proposition}
For any $A, B, C \in \Sigma$, the following results hold true.\\
(1) $\mathcal{M}(A \veebar B) = \mathcal{M}(B \veebar A)$, $\mathcal{M}(A \barwedge B) = \mathcal{M}(B \barwedge A)$;\\
(2) $\mathcal{M}((A \veebar B) \veebar C) = \mathcal{M}(A \veebar (B \veebar C))$, $\mathcal{M}((A \barwedge B) \barwedge C) = \mathcal{M}(A \barwedge (B \barwedge C))$;\\
(3) $\mathcal{M}(A \veebar \bot) = \mathcal{M}(A)$, $\mathcal{M}(A \barwedge \bot) = \mathcal{M}(\bot)$;\\
(4) $\mathcal{M}(A \barwedge \top) = \mathcal{M}(A)$, $\mathcal{M}(A \veebar \top) = \mathcal{M}(\top)$.
\end{proposition}
\begin{proposition}
In a regular vague membership space, for any $A, B, C \in \Sigma$, the following results hold true.\\
(1) $\mathcal{M}(\neg \neg A) = \mathcal{M}(A)$;\\
(2) $\mathcal{M}(\neg (A \veebar B)) = \mathcal{M}((\neg A) \barwedge (\neg B))$, $\mathcal{M}(\neg (A \barwedge B)) = \mathcal{M}((\neg A) \veebar (\neg B))$.
\end{proposition}
\begin{proposition}
In a regular vague membership space, let $A \in \Sigma$. If $\oplus$ and $\otimes$ are complemented with respect to $N$, That is, for any $x \in [0, 1]$, $x \oplus x^{N} = 1$ and $x \otimes x^{N} = 0$, then $\mathcal{M}(\neg A \veebar A) = 1$ and $\mathcal{M}(\neg A \barwedge A) = 0$.
\end{proposition}
\begin{proposition}
In a regular vague membership space, let $N$ be the standard negation and $\oplus = S_{L}, \otimes = T_{L}$. Then for any $A, B \in \Sigma$, satisfies\\
(1) $\mathcal{M}(A \veebar (\neg A \barwedge B)) = \mathcal{M}(A) \vee \mathcal{M}(B)$;\\
(2) $\mathcal{M}(A \barwedge (\neg A \veebar B)) = \mathcal{M}(A) \wedge \mathcal{M}(B)$.
\end{proposition}
\begin{proof}
\begin{eqnarray*}
\mathcal{M}(A \veebar (\neg A \barwedge B)) & = & S_{L}\big(\mathcal{M}(A), \mathcal{M}(\neg A \barwedge B)\big)\\
& = & S_{L}\big(\mathcal{M}(A), T_{L}\big(1 - \mathcal{M}(A), \mathcal{M}(B)\big)\big)\\
& = & \big(\mathcal{M}(A) + ((1- \mathcal{M}(A) + \mathcal{M}(B) - 1) \vee 0)\big) \wedge 1\\
& = & \mathcal{M}(A) \vee \mathcal{M}(B).
\end{eqnarray*}
Hence, (1) holds true. Since
\begin{eqnarray*}
\mathcal{M}(A \barwedge (\neg A \veebar B)) & = & T_{L}\big(\mathcal{M}(A), S_{L}\big(1 - \mathcal{M}(A), \mathcal{M}(B)\big)\big)\\
& = & 1 - S_{L}\big(1 - \mathcal{M}(A), 1 - S_{L}\big(1 - \mathcal{M}(A), \mathcal{M}(B)\big)\big)\\
& = & 1 - S_{L}\big(1 - \mathcal{M}(A), T_{L}\big(\mathcal{M}(A), 1 - \mathcal{M}(B)\big)\big)\\
& = & 1 - \mathcal{M}(\neg A \veebar (A \barwedge \neg B)).
\end{eqnarray*}
Hence, (2) can be obtained from (1) and Proposition 5.3.
\end{proof}
\begin{proposition}
In a regular vague membership space, let $N$ be the standard negation. Then for $p, q \in \Omega$ and $p \neq q$, $\mathcal{M}(p)$ and $\mathcal{M}(q)$ can't all be bigger than 0.5.
\end{proposition}
\begin{proof}
It is a direct conclusion from Axioms III and V.
\end{proof}
\begin{proposition}
Define a binary relation $\equiv$ on $\Sigma$ as follows: for any $A, B \in \Sigma$,
$$A \equiv B\  \text{if and only if}\ \mathcal{M}(A) = \mathcal{M}(B).$$ Then $\equiv$ is an equivalence relation on $\Sigma$. Let $\overline{\Sigma}$ be the set of equivalence classes with respect to $\equiv$, define a binary relation $\precsim$ on $\overline{\Sigma}$ as follows: $$[A] \precsim [B]\  \text{if and only if}\ \mathcal{M}(A) \leqslant \mathcal{M}(B),$$ where $[A]$ denotes the $\equiv$-equivalence class of $A$. Then $(\overline{\Sigma}; \precsim)$ is a chain. If $\oplus = S_{M}$, then $(\overline{\Sigma}; \precsim)$ is a complete chain.
\end{proposition}
\begin{proof}
From the definitions, it is obvious that $\equiv$ is equivalence relation on $\Sigma$, and $\precsim$ is a totally order relation on $\overline{\Sigma}$. Hence, $(\overline{\Sigma}; \precsim)$ is a chain.
\par Let $\oplus = S_{M}$. Next, we need to prove that $(\overline{\Sigma}; \precsim)$ is complete. That is, for every nonempty subset $\{[A_{i}]; i \in I\}$ of $\overline{\Sigma}$ that has an upper bound, has a least upper bound.
\par Since $\veebar_{i \in I}A_{i} \in \Sigma$, it follows from Axiom IV that $[\veebar_{i \in I}A_{i}]$ is an upper bound of $\{[A_{i}]; i \in I\}$. Let $[B]$ be any upper bound of $\{[A_{i}]; i \in I\}$, then $\mathcal{M}(A_{i}) \leqslant \mathcal{M}(B)$ for any $i \in I$. Note that $\mathcal{M}\big(\veebar_{i \in I}A_{i}\big) = \sup\{\mathcal{M}(A_{i}); i \in I\}$ by $\oplus = S_{M}$, thus $\mathcal{M}\big(\veebar_{i \in I}A_{i}\big) \leqslant \mathcal{M}(B)$, that is, $[\veebar_{i \in I}A_{i}] \precsim [B]$. Therefore, $[\veebar_{i \in I}A_{i}]$ is the least upper bound of $\{[A_{i}]; i \in I\}$.
\end{proof}
\begin{proposition}
(The lower limit theorem)Let $[A_{1}] \precsim [A_{2}] \precsim \cdots \precsim [A_{n}] \precsim \cdots$ be an increasing sequence of $\overline{\Sigma}$. Then the limit of the sequence $\{\mathcal{M}(A_{n})\}$ exists. If $\oplus = S_{M}$, then $\lim_{n \rightarrow \infty}\mathcal{M}(A_{n}) = \mathcal{M}\big(\veebar_{n = 1}^{\infty}A_{n}\big)$; otherwise, $\lim_{n \rightarrow \infty}\mathcal{M}(A_{n}) \leqslant \mathcal{M}\big(\veebar_{n = 1}^{\infty}A_{n}\big)$.
\end{proposition}
\begin{proof}
Since $[A_{1}] \precsim [A_{2}] \precsim \cdots \precsim [A_{n}] \precsim \cdots$, it follows from Proposition 5.5 and Axiom I that $$0 \leqslant \mathcal{M}(A_{1}) \leqslant \mathcal{M}(A_{2}) \leqslant \cdots \leqslant \mathcal{M}(A_{n}) \leqslant \cdots \leqslant 1.$$ By the monotone convergence theorem of sequence limit, we have that the limit of the sequence $\{\mathcal{M}(A_{n})\}$ exists and $\lim_{n \rightarrow \infty}\mathcal{M}(A_{n}) = \sup\{\mathcal{M}(A_{i}); i = 1, 2, \cdots\}$. If $\oplus = S_{M}$, then $\mathcal{M}\big(\veebar_{i \in I}A_{i}\big) = \sup\{\mathcal{M}(A_{i}); i = 1, 2, \cdots\}$, and thus $\lim_{n \rightarrow \infty}\mathcal{M}(A_{n}) = \mathcal{M}\big(\veebar_{n = 1}^{\infty}A_{n}\big)$. Otherwise, since $S_{M}$ is the weakest of all t-conorms, thus $$\lim_{n \rightarrow \infty}\mathcal{M}(A_{n}) \leqslant \mathcal{M}\big(\veebar_{n = 1}^{\infty}A_{n}\big).$$
\end{proof}
\par Similarly, it is easy to obtain the following conclusion.
\begin{proposition}
(The upper limit theorem)Let $[A_{1}] \succsim [A_{2}] \succsim \cdots \succsim [A_{n}] \succsim \cdots$ be an decreasing sequence of $\overline{\Sigma}$. Then the limit of the sequence $\{\mathcal{M}(A_{n})\}$ exists. If $\otimes = T_{M}$, then $\lim_{n \rightarrow \infty}\mathcal{M}(A_{n}) = \mathcal{M}\big(\barwedge_{n = 1}^{\infty}A_{n}\big)$; otherwise, $\lim_{n \rightarrow \infty}\mathcal{M}(A_{n}) \geqslant \mathcal{M}\big(\barwedge_{n = 1}^{\infty}A_{n}\big)$.
\end{proposition}
\par In probability theory, the notion of independence is distinguishing and fundamental. However, the case here is different. In the same membership space, we think that any two different vague attributes are related to one another. Hence, the independence of vague attributes can only be defined in two different vague membership spaces. Vague attributes $A$ and $B$ are said to be independent if and only if $A$ and $B$ belong to two different vague spaces with respect to two different attributes, respectively.
\begin{remark}
In fact, the independence in probability theory is also defined in different sample spaces. The only difference is that in the probability theory, the probability measure in the product space of probability spaces $(\Omega_{1}, \Sigma_{1}, P_{1})$ and $(\Omega_{2}, \Sigma_{2}, P_{2})$ is defined as $P_{1} \times P_{2}$; but for vague membership spaces, this way will not work. The membership measure in product vague membership space will be defined in next Section.
\end{remark}
\par Let $(\Omega, \Sigma, \mathcal{M}_{x}, \mathcal{F}, \mathcal{T})$ be a membership space and $A, B \in \Sigma$. Then the conditional membership degree $A$ given $B$ associated with $x$, denoted by $\mathcal{M}_{x}(B \Rightarrow A)$, is defined as $$\mathcal{M}_{x}(B \Rightarrow A) = \sup\{z \in [0, 1] | \mathcal{M}_{x}(B) \otimes z \leqslant \mathcal{M}_{x}(A)\}.$$
The value $\mathcal{M}_{x}(B \Rightarrow A)$ is the measure of relationship between vague attributes $B$ and $A$  associated with $x$, and show the degree of vague attribute $B$ influences vague attribute $A$ in the membership space $(\Omega, \Sigma, \mathcal{M}_{x}, \mathcal{F}, \mathcal{T})$. Furthermore, we can define the conditional membership degree $A$ given $B$, denoted by $A | B$, in vague space $(\Omega, \Sigma, \mathcal{F})$ as follows:
$$A | B = \min_{x \in U}\{\mathcal{M}_{x}(B \Rightarrow A)\}.$$
\section{Product vague membership space}
\par Thus far, we have established the membership degree theory for one dimensional vague space by axiomatical approach. More general, for multidimensional vague attributes or multidimensional vague predicates, we need to establish the product vague membership space, we will take the Cartesian product $[0, 1]^{n}$ as the range of product vague membership measure.
\par As mentioned early, a multidimensional attribute can be determined usually by several one dimensional attributes, and can be regarded as the Cartesian product of these one-dimensional attributes. Therefore, we define product vague membership space as the Cartesian product of several one dimensional vague membership spaces.
\begin{definition}
Let $(\Omega_{i}, \Sigma_{i}, \mathcal{F}_{i}), i = 1, 2, \cdots, n$, be vague spaces.\\
(1) $\Omega_{1} \times \cdots \times \Omega_{n} \equiv \{(\omega_{1}, \omega_{2}, \cdots, \omega_{n}) | \omega_{i} \in \Omega_{i}, i = 1, 2, \cdots, n\}$, the set of all ordered arrays, is called product elementary vague attribute set of $\Omega_{1}$, $\cdots$, $\Omega_{n}$.\\
(2) The product vague attribute set of $\Sigma_{1}$, $\cdots$, $\Sigma_{n}$ on $\Omega_{1} \times \cdots \times \Omega_{n}$, denoted by $\Sigma_{1} \times \cdots \times \Sigma_{n}$, is the set of all ordered arrays, i.e., $$\Sigma_{1} \times \cdots \times \Sigma_{n} \equiv \{(A_{1}, \cdots, A_{n})| A_{i} \in \Sigma_{n}, i = 1, 2, \cdots, n\}.$$
Every element of $\Sigma_{1} \times \cdots \times \Sigma_{n}$ is called a $n$-dimensional vague attributes.
\end{definition}
\par Given a concept $C$, consider the partition of the extension (which is usual a crisp (or classical) set or class) of $C$ according to a group of one dimensional attributes $\varphi_{C}^{1}, \cdots, \varphi_{C}^{n}$, $n \in \mathbb{N}^{+}$. Let $(\Omega_{\varphi_{C}^{i}}, \Sigma_{\varphi_{C}^{i}}, \mathcal{M}_{x_{i}}^{i}, \mathcal{F}_{i}, \mathcal{T}_{i})$ be the vague membership space with respect to $\varphi_{C}^{i}$ and associated with $x_{i} \in U_{\varphi_{C}^{i}}$, then the vague membership space with respect to this group attributes $\Psi_{C} = \{\varphi_{C}^{1}, \cdots, \varphi_{C}^{n}\}$ and associated with $\overline{x} = (x_{1}, \cdots, x_{n})$ is such a quintuple $(\Omega_{\Psi_{C}}, \Sigma_{\Psi_{C}}, \mathcal{M}_{\overline{x}}, \mathcal{F}, \mathcal{T})$ which is called product vague membership space of $(\Omega_{\varphi_{C}^{i}}, \Sigma_{\varphi_{C}^{i}}, \mathcal{M}_{x_{i}}^{i}, \mathcal{F}_{i}, \mathcal{T}_{i})$, where $\Omega_{\Psi_{C}} = \Omega_{\varphi_{C}^{1}} \times \Omega_{\varphi_{C}^{2}} \times \cdots \times \Omega_{\varphi_{C}^{n}}$ and $\Sigma_{\Psi_{C}} = \Sigma_{\varphi_{C}^{1}} \times \cdots \times \Sigma_{\varphi_{C}^{n}}$. Let $\mathcal{F}_{i} = \{\bot_{i}, \top_{i}, \neg_{i}, \veebar_{i}, \barwedge_{i}\}$, $\mathcal{T}_{i} = \{N_{i}, \oplus_{i}, \otimes_{i}\}$, $i = 1, 2, \cdots, n$, then $\mathcal{F} = \mathcal{F}_{1} \times \cdots \times \mathcal{F}_{n}$, and $\mathcal{T} = \mathcal{T}_{1} \times \cdots \times \mathcal{T}_{n}$. $\mathcal{M}_{\overline{x}}$ is a real vector valued function from $\Sigma_{\Psi_{C}}$ to $[0, 1]^{n}$ with respect to $\overline{x}$ such that $$\mathcal{M}_{\overline{x}}((A_{1}, \cdots, A_{n})) = (\mathcal{M}_{x_{1}}^{1}(A_{1}), \cdots, \mathcal{M}_{x_{n}}^{n}(A_{n}))$$ for any $(A_{1}, \cdots, A_{n}) \in \Sigma_{\Psi_{C}}$. Hence,\\
(1) For any $(A_{1}, \cdots, A_{n}) \in \Sigma_{\Psi_{C}}$, $\textbf{0} \leqslant \mathcal{M}_{\overline{x}}((A_{1}, \cdots, A_{n})) \leqslant \textbf{1}$, where $\textbf{0} = (\overbrace{0, \cdots, 0}^{n})$ and $\textbf{1} = (\overbrace{1, \cdots, 1}^{n})$.\\
(2) Let $\bot = (\bot_{1}, \cdots, \bot_{n})$ and $\top = (\top_{1}, \cdots, \top_{n})$. Then $\mathcal{M}_{\overline{x}}(\bot) = \textbf{0}$, $\mathcal{M}_{\overline{x}}(\top) = \textbf{1}$.\\
(3) For any $(A_{1}, \cdots, A_{n}) \in \Sigma_{\Psi_{C}}$, $$\mathcal{M}_{\overline{x}}((\neg_{1} A_{1}, \cdots, \neg_{n} A_{n})) \leqslant \big((\mathcal{M}^{1}_{x_{1}}(A_{1}))^{N_{1}}, \cdots, (\mathcal{M}^{n}_{x_{n}}(A_{n}))^{N_{n}}\big).$$
(4) For any $(A_{1i}, \cdots, A_{ni}),  i = 1, 2, \cdots, \in \Sigma_{\Psi_{C}}$,
\begin{eqnarray*}
& & \mathcal{M}_{\overline{x}}\big((A_{11} \veebar_{1} A_{12} \veebar_{1} \cdots, \cdots, A_{n1} \veebar_{n} A_{n2} \veebar_{n} \cdots)\big) \\
& = & \big(\mathcal{M}^{1}_{x_{1}}(A_{11}) \oplus_{1} \mathcal{M}^{1}_{x_{1}}(A_{12}) \oplus_{1} \cdots, \cdots, \mathcal{M}^{n}_{x_{n}}(A_{n1}) \oplus_{n} \mathcal{M}^{n}_{x_{n}}(A_{n2}) \oplus_{n} \cdots \big).
\end{eqnarray*}
(5) For any $(A_{1i}, \cdots, A_{ni}),  i = 1, 2, \cdots, \in \Sigma_{\Psi_{C}}$,
\begin{eqnarray*}
& & \mathcal{M}_{\overline{x}}(A_{11} \barwedge_{1} A_{12} \barwedge_{1} \cdots, \cdots, A_{n1} \barwedge_{n} A_{n2} \barwedge_{n} \cdots) \\
& = & \big(\mathcal{M}^{1}_{x_{1}}(A_{11}) \otimes_{1} \mathcal{M}^{1}_{x_{1}}(A_{12}) \otimes_{1} \cdots, \cdots, \mathcal{M}^{n}_{x_{n}}(A_{n1}) \otimes_{n} \mathcal{M}^{n}_{x_{n}}(A_{n2}) \otimes_{n} \cdots \big).
\end{eqnarray*}
\par The $n$-dimensional product vague membership space $(\Omega_{\Psi_{C}}, \Sigma_{\Psi_{C}}, \mathcal{M}_{\overline{x}}, \mathcal{F}, \mathcal{T})$ is said to be regular if we take the equality in (3), and is said to be normal if there exists $(p_{1}, \cdots, p_{n}) \in \Omega_{\Psi_{C}}$ such that $$\mathcal{M}_{\overline{x}}((p_{1}, \cdots, p_{n})) = \textbf{1}.$$ In other words, $(\Omega_{\Psi_{C}}, \Sigma_{\Psi_{C}}, \mathcal{M}_{\overline{x}}, \mathcal{F}, \mathcal{T})$ is normal if and only if every component vague membership space $(\Omega_{\varphi_{C}^{i}}, \Sigma_{\varphi_{C}^{i}}, \mathcal{M}_{x_{i}}^{i}, \mathcal{F}_{i}, \mathcal{T}_{i})$ is normal, $i = 1, 2, \cdots, n$. Then, $\overline{x}$ is said to be crisp in product vague membership space $(\Omega_{\Psi_{C}}, \Sigma_{\Psi_{C}}, \mathcal{M}_{\overline{x}}, \mathcal{F}, \mathcal{T})$. Similar to the case of one dimensional, any partition in the classical sense of Cartesian product $X^{n}$ with respect to any $\overline{x} \in X^{n}$ is a $n$-dimensional normal product vague membership space.
\par In what follows, for the sake of convenience, we denote a $n$-dimensional product vague membership space $(\Omega_{\Psi_{C}}, \Sigma_{\Psi_{C}}, \mathcal{M}_{\overline{x}}, \mathcal{F}, \mathcal{T})$ by $(\Omega^{n}, \Sigma^{n}, \mathcal{M}^{n}, \mathcal{F}^{n}, \mathcal{T}^{n})$ when it does not involve any concrete problems. Moreover, the operations in $\mathcal{F}^{n}$ and the corresponding operations $\mathcal{T}^{n}$ are defined pointwisely on $\Sigma^{n}$ and on $[0, 1]^{n}$, respectively.
\begin{proposition}
Let $(\Omega^{n}, \Sigma^{n}, \mathcal{M}^{n}, \mathcal{F}^{n}, \mathcal{T}^{n})$ be a $n$-dimensional product vague membership space. Define a binary relation $\equiv$ on $\Sigma^{n}$ as follows: for any $A = (A_{1}, \cdots, A_{n}), B = (B_{1}, \cdots, B_{n}) \in \Sigma^{n}$,
$$A \equiv B\  \text{if and only if}\ \mathcal{M}_{i}(A_{i}) = \mathcal{M}_{i}(B_{i}), i = 1, 2, \cdots, n.$$ Then $\equiv$ is an equivalence relation on $\Sigma^{n}$. Let $\overline{\Sigma^{n}}$ be the set of equivalence classes with respect to $\equiv$, define a binary relation $\precsim$ on $\overline{\Sigma^{n}}$ as follows: $$[A] \precsim [B]\  \text{if and only if}\ \mathcal{M}_{i}(A_{i}) \leqslant \mathcal{M}_{i}(B_{i}), i = 1, 2, \cdots, n.$$ where $[A]$ denotes the $\equiv$-equivalence class of $A$. If $\oplus = (\oplus_{1}, \cdots, \oplus_{n}) = (S_{M}, \cdots, S_{M})$, then $(\overline{\Sigma^{n}}; \precsim)$ is a complete lattice.
\end{proposition}
\begin{proof}
Similar to Proposition 5.6, it is easy to prove that $\equiv$ is equivalence relation on $\Sigma^{n}$, and $\precsim$ is a partial order relation on $\overline{\Sigma^{n}}$.
\par For any $A = (A_{1}, \cdots, A_{n}), B = (B_{1}, \cdots, B_{n}) \in \Sigma^{n}$, let $$C = A \veebar B = (A_{1} \veebar_{1} B_{1}, \cdots, A_{n} \veebar_{n} B_{n}).$$ Since $\oplus = (\oplus_{1}, \cdots, \oplus_{n}) = (S_{M}, \cdots, S_{M})$, it follows that $$\mathcal{M}_{i}(A_{i} \veebar_{i} B_{i}) = \mathcal{M}_{i}(A_{i}) \oplus_{i} \mathcal{M}_{i}(B_{i}) = \max\{\mathcal{M}_{i}(A_{i}), \mathcal{M}_{i}(B_{i})\}, i = 1, 2, \cdots, n,$$ thus $[A] \precsim [C]$ and $[B] \precsim [C]$, $[C]$ is an upper bound of $[A]$ and $[B]$ with respect to $\precsim$. If $D = (D_{1}, \cdots, D_{n})$, $[D]$ is also an upper bound of $[A]$ and $[B]$, then $\mathcal{M}_{i}(A_{i}) \leqslant \mathcal{M}_{i}(D_{i})$ and $\mathcal{M}_{i}(B_{i}) \leqslant \mathcal{M}_{i}(D_{i})$, $i = 1, 2, \cdots, n$. Hence, $\max\{\mathcal{M}_{i}(A_{i}), \mathcal{M}_{i}(B_{i})\} = \mathcal{M}_{i}(C_{i}) \leqslant \mathcal{M}_{i}(D_{i})$, $i = 1, 2, \cdots, n$, that is, $[C] \precsim [D]$. Therefore, $[C]$ is the supremum of $[A]$ and $[B]$. Similarly, we can prove that $[A \barwedge B]$ is the infimum of $[A]$ and $[B]$. So far, we have proved that $(\overline{\Sigma^{n}}; \precsim)$ is a lattice.
\par Note that for any $\{[A^{i}]; i \in I\} \subset \overline{\Sigma^{n}}$, $\veebar_{i \in I}[A^{i}], \barwedge_{i \in I}[A^{i}] \in \overline{\Sigma^{n}}$, the completeness can be proved similarly by Axiom IV.
\end{proof}
\par Similar to Proposition 5.7 and 5.8, we can prove the following conclusions.
\begin{proposition}
(The lower limit theorem)Let $[A^{1}] \precsim [A^{2}] \precsim \cdots \precsim [A^{m}] \precsim \cdots$ be an increasing sequence of $\overline{\Sigma^{n}}$. Then the limit of the vector valued sequence $\{\mathcal{M}(A^{m})\}$ exists. If $\oplus = (\oplus_{1}, \cdots, \oplus_{n}) = (S_{M}, \cdots, S_{M})$, then $$\lim_{m \rightarrow \infty}\mathcal{M}(A^{m}) = \Big(\mathcal{M}_{1}\big((\veebar_{1})_{m = 1}^{\infty}A^{m}_{1}\big), \cdots, \mathcal{M}_{n}\big((\veebar_{n})_{m = 1}^{\infty}A^{m}_{n}\big)\Big).$$ Otherwise, $\lim_{m \rightarrow \infty}\mathcal{M}(A^{m}) \leqslant \Big(\mathcal{M}_{1}\big((\veebar_{1})_{m = 1}^{\infty}A^{m}_{1}\big), \cdots, \mathcal{M}_{n}\big((\veebar_{n})_{m = 1}^{\infty}A^{m}_{n}\big)\Big)$.
\end{proposition}
\begin{proposition}
(The upper limit theorem)Let $[A^{1}] \succsim [A^{2}] \succsim \cdots \succsim [A^{m}] \succsim \cdots$ be an decreasing sequence of $\overline{\Sigma^{n}}$. Then the limit of the vector valued sequence $\{\mathcal{M}(A^{n})\}$ exists. If $\otimes = (\otimes_{1}, \cdots, \otimes_{n}) = (T_{M}, \cdots, T_{M})$, then
$$\lim_{m \rightarrow \infty}\mathcal{M}(A^{m}) = \Big(\mathcal{M}_{1}\big((\barwedge_{1})_{m = 1}^{\infty}A^{m}_{1}\big), \cdots, \mathcal{M}_{n}\big((\barwedge_{n})_{m = 1}^{\infty}A^{m}_{n}\big)\Big).$$ Otherwise, $\lim_{m \rightarrow \infty}\mathcal{M}(A^{m}) \geqslant \Big(\mathcal{M}_{1}\big((\barwedge_{1})_{m = 1}^{\infty}A^{m}_{1}\big), \cdots, \mathcal{M}_{n}\big((\barwedge_{n})_{m = 1}^{\infty}A^{m}_{n}\big)\Big)$.
\end{proposition}
\section{Vague variables and vague vectors}
In this section, we introduce a very important notion in vague membership space, vague variable, which is similar to the notion of random variable in the probability theory. Vague variables will be an important tool for studying phenomenons of vagueness.
\begin{definition}
Let $(\Omega, \Sigma, \mathcal{M}, \mathcal{F}, \mathcal{T})$ be a vague membership space. The real valued function $X: \Omega \cup \{\top, \bot\} \rightarrow \mathbb{R}$ is called a vague variable on $(\Omega, \Sigma, \mathcal{M}, \mathcal{F}, \mathcal{T})$ if $X(\top) = + \infty$ and $X(\bot) = - \infty$.
\end{definition}
\begin{example}
Let $\Omega = \{Young, Medium, Old, Senior\}$. Define the function $X = \text{"Older man"}$ as follows: $X(\top) = + \infty$, $X(\bot) = - \infty$,
$$X(Young) = 0,\ X(Medium) = 1,\ X(Old) = 1, \text{  and } X(Senior) = 1,$$
then $X$ is a vague variable.
\end{example}
\begin{example}
Let $\Omega = \{\text{Cold, Cool, Normal temperature, Warm, Hot, Tropical}\}$. Define the mapping $X = \text{"Acceptable temperature"}$ by $X(\top) = + \infty$, $X(\bot) = - \infty$, $X(Cold) = 0$, $X(Cool) = 1$, $X(\text{Normal temperature}) = 1$, $X(Warm) = 1$, $X(Hot) = 0$, $X(Tropical) = 0$. Then $X$ is a vague variable.
\end{example}
\begin{definition}
Let $X$ be a vague variable on $(\Omega, \Sigma, \mathcal{M}, \mathcal{F}, \mathcal{T})$. Let $$F_{X}(x) \equiv \bigoplus_{p \in \{X \leqslant x\}}\mathcal{M}(p), x \in \mathbb{R},$$ where $\{X \leqslant x\} \triangleq \{p \in \Omega \cup \{\top, \bot\}: X(p) \leqslant x\}$. Then $F_{X}(\cdot)$ is called the membership degree cumulative distribution function of $X$ on $(\Omega, \Sigma, \mathcal{M}, \mathcal{F}, \mathcal{T})$.
\end{definition}
\par In what follows, we always write $\{X \leqslant x\}$ instead of $\{p \in \Omega \cup \{\top, \bot\}: X(p) \leqslant x\}$ when it leads to no confusion.
\begin{theorem}
Let $X$ be a vague variable and $F_{X}(\cdot)$ the membership degree cumulative distribution function of $X$ on $(\Omega, \Sigma, \mathcal{M}, \mathcal{F}, \mathcal{T})$. The following properties are satisfied:\\
F1. Nondecreasing properties: for any $x_{1}, x_{2} \in \mathbb{R}$, if $x_{1} \leqslant x_{2}$, then $F_{X}(x_{1}) \leqslant F_{X}(x_{2})$;\\
F2. Left continuity: if $\oplus$ is continuous, then for any $x_{0} \in \mathbb{R}$, $\lim_{t \rightarrow x_{0}^{-}}F_{X}(t) = F_{X}(x_{0})$;\\
F3. $\lim_{t \rightarrow -\infty}F_{X}(t) = 0, \lim_{t \rightarrow \infty}F_{X}(t) = 1$.
\end{theorem}
\begin{proof}
Since for $x_{1} \leqslant x_{2}$ one has $\{X \leqslant x_{1}\} \subseteq \{X \leqslant x_{2}\}$, F1 immediately follows from Definition 4.2 and 7.2.
\par Let $\{x_{n}\}$ be an increasing sequence with $x_{n} \rightarrow x_{0}$, $A = \{X \leqslant x_{0}\}$, $A_{n} = \{X \leqslant x_{n}\}$. Then the sequence of sets $A_{n}$ also increases, and $\bigcup A_{n} = A$. Therefore, $$\bigoplus_{p \in \{X \leqslant x_{n}\}}\mathcal{M}(p) \rightarrow \bigoplus_{p \in \{X \leqslant x_{0}\}}\mathcal{M}(p).$$ This means that $\lim_{t \rightarrow x_{0}^{-}}F_{X}(t) = F_{X}(x_{0})$. F2 holds.
\par To prove F3, consider two number sequences $\{x_{n}\}$ and $\{y_{n}\}$ such that $\{x_{n}\}$ is decreasing and $x_{n} \rightarrow - \infty$, while $\{y_{n}\}$ is increasing and $y_{n} \rightarrow \infty$. Put $A_{n} = \{X \leqslant x_{n}\}$, $B_{n} = \{X \leqslant y_{n}\}$. Since $x_{n}$ tends monotonically to $- \infty$, the sequence of sets $A_{n}$ decreases monotonically to $\bigcap A_{n} = \{\bot\}$. Therefore, $\bigoplus_{p \in \{X \leqslant x_{n}\}}\mathcal{M}(p) \rightarrow \mathcal{M}(\bot) = 0$. That is, $\lim_{t \rightarrow -\infty}F_{X}(t) = 0$. Since the sequence $\{y_{n}\}$ tends monotonically to $\infty$, the sequence of sets $B_{n}$ increases to $\bigcup B_{n} = \Omega \cup \{\top, \bot\}$. This implies, as above, that $\lim_{n \rightarrow \infty}F_{X}(y_{n}) = 1$, $\lim_{x \rightarrow \infty}F_{X}(x) = 1$.
\end{proof}
\begin{definition}
(1) A vague variable $X$ on $(\Omega, \Sigma, \mathcal{M}, \mathcal{F}, \mathcal{T})$ is called discrete if $\Omega$ is a finite or countable set.\\
(2) A vague variable $X$ on $(\Omega, \Sigma, \mathcal{M}, \mathcal{F}, \mathcal{T})$ is called continuous if $\Omega$ is a continuous set.
\end{definition}
\begin{definition}
Let $(\Omega, \Sigma, \mathcal{M}, \mathcal{F}, \mathcal{T})$ be a vague membership space, $k \in \mathbb{N}^{+}$. The vector valued function $X = (X_{1}, \cdots, X_{k}): \Omega \cup \{\top, \bot\} \rightarrow \mathbb{R}^{k}$ is called a (k-dimensional) vague vector on $(\Omega, \Sigma, \mathcal{M}, \mathcal{F}, \mathcal{T})$ if $X_{i}$ is a vague variable on $(\Omega, \Sigma, \mathcal{M}, \mathcal{F}, \mathcal{T})$ for any $i = 1, 2, \cdots, k$.
\end{definition}
\par Let $X = (X_{1}, \cdots, X_{k})$ be a vague vector with components $X_{i}$, $i = 1, 2, \cdots, k$. Then each $X_{i}$ is a vague variable on $(\Omega, \Sigma, \mathcal{M}, \mathcal{F}, \mathcal{T})$. Conversely, if for $1 \leqslant i \leqslant k$, $X_{i}$ is a vague variable on $(\Omega, \Sigma, \mathcal{M}, \mathcal{F}, \mathcal{T})$, then $X = (X_{1}, \cdots, X_{k})$ is a vague vector.
\begin{definition}
Let $X$ be a $k$-dimensional vague vector on $(\Omega, \Sigma, \mathcal{M}, \mathcal{F}, \mathcal{T})$ for some $k \in \mathbb{N}$. Let $$F_{X}(\overline{x}) \equiv \bigoplus_{p \in \{X \leqslant \overline{x}\}}\mathcal{M}(p)$$ for $\overline{x} = (x_{1}, x_{2}, \cdots, x_{k}) \in \mathbb{R}^{k}$, where $$\{X \leqslant \overline{x}\} \triangleq \{p \in \Omega \cup \{\top, \bot\}: X_{1}(p) \leqslant x_{1}, X_{2}(p) \leqslant x_{2}, \cdots, X_{k}(p) \leqslant x_{k}\}.$$ Then $F_{X}(\cdot)$ is called the membership degree joint cumulative distribution function of the vague vector $X$ on $(\Omega, \Sigma, \mathcal{M}, \mathcal{F}, \mathcal{T})$.
\end{definition}
\begin{definition}
Let $X = (X_{1}, \cdots, X_{k})$ be a vague vector on $(\Omega, \Sigma, \mathcal{M}, \mathcal{F}, \mathcal{T})$. Then,  for each $i = 1, \cdots, k,$, the membership degree cumulative distribution function $F_{X_{i}}$ of the vague variable $X_{i}$ is called the marginal membership degree cumulative distribution function of $X_{i}$.
\end{definition}
\par It is clear that the membership degree joint cumulative distribution $F_{X}$ of $X$ determines the marginal membership degree cumulative distribution $F_{X_{i}}$ of $X_{i}$ for all $i = 1, 2, \cdots, k$. However, the marginal membership degree cumulative distributions $\{F_{X_{i}}: i = 1, 2, \cdots, k\}$ do not uniquely determine the joint membership degree cumulative distribution $F_{X}$, without additional conditions.
\par Next, we define the notion of balanced (acceptable) value of a vague variable, which is similar to expected value of a random variable in probability theory.
\begin{definition}
Let $X$ be a vague variable on $(\Omega, \Sigma, \mathcal{M}, \mathcal{F}, \mathcal{T})$. The balanced (or acceptable) value of $X$, denoted by $B(X)$, is defined as $$B(X) = \frac{\int_{\Omega}Xd\mathcal{M}}{\int_{\Omega}d\mathcal{M}},$$ provided the integral is well defined. That is, at least one of the two quantities $\int X^{+}d\mathcal{M}$ and $\int X^{+}d\mathcal{M}$ is finite.
\end{definition}
\begin{remark}
If $X$ is a discrete vague variable on $(\Omega, \Sigma, \mathcal{M}, \mathcal{F}, \mathcal{T})$, then $$B(X) = \frac{\sum_{p \in \Omega}X(p)\mathcal{M}(p)}{\sum_{p \in \Omega}\mathcal{M}(p)}.$$
\end{remark}
\section{Vague partition and fuzzy sets}
\par In the preceding sections, we presented the axioms for membership degrees, and defined some fundamental notions, such as vague variable, etc., then you may wonder what the relationship between these elements and Zadeh's fuzzy set is.
\par In this section, we will introduce the notion of vague partition, and show that this relationship between vague variables and vague partitions is similar to that of random variables and stochastic processes in probability theory, and that fuzzy sets can be generated by vague partitions.
\begin{definition}
Let $U$ be the domain of discourse and $\Omega$ an elementary vague attribute set. A vague partition of $U$ is a family $\{X(p, x): p \in \Omega, x \in U\}$ of vague variables defined on $\Omega$. For every $x \in U$, $X(p, x)$ is a vague variable in the vague membership space $(\Omega, \Sigma, \mathcal{M}_{x}, \mathcal{F}, \mathcal{T})$. And for every $p \in \Omega$, $X(p, x)$ is a vague class or a fuzzy subset in $U$ with respect to $p$.
\end{definition}
\par The vague partition is said to be regular if in Definition 8.1, for every $x \in U$, the vague membership space $(\Omega, \Sigma, \mathcal{M}_{x}, \mathcal{F}, \mathcal{T})$ is regular.
\begin{example}
Let $C$ be the concept "Man", $\varphi_{C}$ the attribute "Age" and $U = [0, 200]$. The elementary vague attribute set $\Omega_{\varphi_{C}}$ according to $\varphi_{C}$ is defined as follows: $$\Omega_{\varphi_{C}} \triangleq \{\text{Childhood, Juvenile, Youth, Maturity, Midlife, Elder, Senectitude}\}.$$ Define the function $X = "\text{Young adults}"$ as follows: $X(\top) = + \infty, X(\bot) = - \infty$, $$X(\text{Childhood}) = X(\text{Juvenile}) = X(\text{Elder}) = X(\text{Senectitude}) = 0,$$  $X(\text{Youth}) = X(\text{Maturity}) = 1$ and $X(\text{Midlife}) = \frac{1}{2}$, then for any $x \in U$, $X$ is a vague variable in the vague membership space $(\Omega, \Sigma, \mathcal{M}_{x}, \mathcal{F}, \mathcal{T})$, hence $\{X(p, x) : p \in \Omega_{\varphi_{C}}, x \in U\}$ is a regular vague partition of $U$. $X(Youth, x)$ is a fuzzy subset in $U$ with respect to the attribute "Youth", its membership function $\mu_{Youth}$ is defined as $\mu_{Youth}(x) = \mathcal{M}_{x}(Youth)$ for any $x \in U$.
\end{example}
\begin{remark}
In the definition above, it is worth noting that the notion of fuzzy set introduced by Zadeh in 1965 can be redefined by the notion of vague partition. In fact, a vague partition will generate several fuzzy sets. For a fixed vague attribute $p$, the vague partition can generate a fuzzy set with respect to $p$.
\par Moreover, when $\Omega$ is a set of multidimensional vague attributes, we take the Cartesian product $[0,1]^{n}$ as the range of product vague membership measure, then the vague partition $\{X(p, x): p \in \Omega, x \in U\}$ of $U$ can generate $L$-fuzzy sets.
\end{remark}
\begin{remark}
In addition, various extensions of fuzzy sets, such as, intuitionistic fuzzy set and vague sets, can also be incorporated into my framework according to Axiom III. In fact, intuitionistic fuzzy set \cite{Atan86} can only be regarded as a particular case in our axiomatic system, where the set of elementary vague attributes consists only of one element.
\end{remark}
\par According to definition 8.1 and the axiomatic system in section 4, the membership function of a fuzzy set generated by some vague partition should satisfy these axioms, hence it will be restricted by the membership functions of other fuzzy sets generated by the vague partition. In such a way, the obtained membership function should be more objective.
\par Moreover, it is easy to find that the shape of the membership function of a fuzzy set is determined not only by these axioms governing membership degrees and their interconnections, but also by the set of elementary vague attributes in the vague space. As we mentioned earlier, the set of elementary vague attributes is changeable.
\par It also need to note that the intersection and union of two fuzzy sets can be defined only in the same vague membership space, fuzzy sets in two different vague membership spaces can not define the operations of intersection and union directly, but they can be considered in product vague membership space.
\par Maybe you have found that we also can consider the time factor into the vague partition just like the notion of stochastic process in probability theory, then we can obtain the notion of vague partition process, which can be used to model the vagueness of in different states, or to model the change of vague phenomena over time. For example, the growth process of a man, the changing process of temperature over time, etc.
\begin{definition}
Let $U$ be the domain of discourse and $\Omega$ an elementary vague attribute set. A vague partition process of $U$ is a family $\{X(p, x, t): p \in \Omega, x \in U, t \in T\}$ of vague variables defined on $\Omega$, where $T$ is the index set. For every $t \in T$, $X(p, x, t)$ is a vague partition of $U$. For every $p \in \Omega$, $X(p, x, t)$ is called an object trajectory of the vague partition process. And for every $x \in U$, $X(p, x, t)$ is called an attribute trajectory of the vague partition process.
\end{definition}
\begin{remark}
In a vague partition process, the parameter $t$ can also be regarded as the possible world (or the existing state) of the corresponding vague attributes, or the cognitive level about these vague attributes. As an example, we consider the vague predicate "\textit{is tall}". Obviously, its meaning is different in different background, e.g., in Venezuela and in Netherlands.
\end{remark}
\par The vague partition process is said to be regular if in Definition 8.2, for every $x \in U$, the vague membership space $(\Omega, \Sigma, \mathcal{M}_{x}, \mathcal{F}, \mathcal{T})$ is regular.
\par A vague partition process is called discrete if $T$ is a finite or countable set, is called continuous if $T$ is a continuous set.

\section{Conclusion}
\par In 1930s, Kolmogorov \cite{Kolmogorov1933} established the rigorous foundations of theory of probability by an axiomatical mathematical formalization, which has made theory of probability an acknowledged and independent mathematical branch. In probability theory, probability is a numerical measure of the likeliness that an event will occur. Probability theory focus more on the relationship among random events in a sample space, rather than the accuracy of the probability values (numbers between 0 and 1), the probabilities of these random events are required to satisfy all of the probability axioms which reflect the intuitive characterization of random phenomenon. This means that probability theory characterizes random phenomena from a global and overall point of view (consider random events in a probability space), not just focuses on a point (a single random event).
\par Although vagueness is different from randomness in their origins as we have mentioned in Section 1, but they are similar in other ways, especially from the perspective of epistemology. The mathematical model treating vagueness aims at modelling vague phenomena approximately in terms of a kind of numerical measure, and it's not our aim to represent vague phenomena accurately. So in this sense, probability theory and vague membership theory are similar. Hence, the successful approaches establishing probability theory can inspire us to establish a suitable vague membership degree theory. It's just with these thoughts in mind that in this paper, we established the axiomatical foundation for the theory of membership degrees which deals with vague phenomenon by many-valued methods.
\par However, we should also be aware that Lebesgue's theories of measure and integration (see \cite{Cohn2013}), which have been used to establish theory of probability, are not suitable to establish theory of membership degree for vagueness. For example, if we consider the conjunction of two vague attributes (such as, \textit{Red} and \textit{Orange}) as the intersection of two sets ($\{Red\}$ and $\{Orange\}$), then this will not work, that is because the intersection is an empty set, but "\textit{Red} and \textit{Orange}" is obvious a new vague attribute. Based on this consideration, in this paper, we employed the methods in universal algebra to establish the basic concepts and principles.
\par The axiomatical model proposed in this paper will help to explain and to standardize the use of operations and extensions of fuzzy sets, to avoid any misuse of fuzzy set theory and to make the methods based fuzzy set theory more scientific and objective. Meanwhile, based on the new axiomatical model, some new theories and methods processing vague phenomena can also be further developed. Meanwhile, this paper also explained the relation between natural language and fuzzy sets.
\par The thesis defended in this paper is that the precondition to model and characterize vague phenomena is to know what attributes bring the vagueness, one dimensional or multidimensional. The difference among vague attributes is the key point to recognize and model vague phenomena, membership degrees should be considered in a vague membership space. In such a way, the notion of fuzzy set has been incorporated into the axiomatic system of membership degrees, so we have established the axiomatical foundation of fuzzy sets.
\par We hope that the work in this paper should provide with an axiomatical mathematical model for dealing with vague phenomena from a global point of view. It's also our hope that this work can serve as a tool to expound those long-standing controversies and divergences in fuzzy set theory and its applications. We think that a good formalized theory or method treating vagueness should have the resources to accommodate all the different types of vague phenomena, and its intuitive meaning is clear.
\section{Acknowledgements}
I would like to express my warm thanks to Prof. Y. Xu, P. Eklund and D.W. Pei for valuable discussions on some of the problems considered here.
\par The work was partially supported by the National Natural Science Foundation of
China (Grant No. 61100046, 61175055, 61305074), the Application Fundamental Research Plan Project of Sichuan Province (Grant No. 2011JY0092) and the Fundamental Research Funds for the Central Universities of China (Grant No.2682014ZT28).





\bibliographystyle{model1b-num-names}
\bibliography{<your-bib-database>}



\end{document}